\RequirePackage{lineno}

    \documentclass[12pt]{article}
\usepackage{amssymb,amsmath,amscd,amsthm}

\newtheorem{theorem}{Theorem}[section]

\newtheorem{lemma}[theorem]{Lemma}

\date{}

\begin{document}
\title{Population dynamics with moderate tails of the underlying random walk }

\author{ S. Molchanov \footnote{Dept of Mathematics and Statistics, UNC at Charlotte, NC 28223 and National Research Univ., Higher School of Economics, Russian Federation, smolchan@uncc.edu. The work was partially supported  by the NSF grant DMS-1714402 and by the Russian Science Foundation, project \# 17-11-01098.}, B. Vainberg \footnote{Dept of Mathematics and Statistics, UNC at Charlotte, NC 28223, brvainbe@uncc.edu. The work was partially supported  by the NSF grant DMS-1714402 and the Simons Foundation grant 527180; corresponding author.}}

\maketitle

\begin{abstract}
Symmetric random walks in $R^d$ and $Z^d$ are considered. It is assumed that the jump distribution density has moderate tails, i.e., several density moments are finite, including the second one. The global (for all $x$ and $t$) asymptotic behavior at infinity of the transition probability (fundamental solution of the corresponding parabolic convolution operator) is found. Front propagation of ecological waves in the corresponding population dynamics models is described.
\end{abstract}



\section{Introduction}
The main object in the population dynamics is the particle field, i.e., the integer valued random measure $n(t,\Gamma)$ defining the number of particles at time $t \geq 0$ in the set $\Gamma\subset R^d$. The evolution of this measure in time can be described as follows.

At the initial moment, one usually has either finitely many particles or an infinite homogeneous in space initial field $n(0,\Gamma)$, e.g.,  the Poisson field with density $\lambda>0$. Each initial particle generates its own subpopulation, independently of others. Starting from the initial location $x\in R^d$, the particle performs a random walk with the generator
\begin{equation}\label{HHHa}
\mathcal L\psi(x)=\chi\int_{ R^d}(\psi(x+z)-\psi(x))a(z)dz, \quad x\in R^d.
\end{equation}
The particle stays in $x$  for a random time $\tau$ with the exponential distribution Exp$(\kappa)$ and then jumps from $x$ to $x+z$ with the distribution density
$a(z),z\in R^d.$ Of course,
\begin{equation}\label{aone1}
a\geq 0, \quad \int_{R^d}a(z)dz=1.
\end{equation}
We will always assume the symmetry of the walk:
\begin{equation}\label{aone1sym}
a(z)=a(-z).
\end{equation}

The evolution of the particles includes also the birth-death processes. In the simplest scenario (no aging phenomena, etc),  during the time interval $(t,t+dt)$, each particle either dies with the probability $\mu dt$  ($\mu$ is the mortality rate)
or splits into two particles with the probability $\beta dt$ ($\beta$ is the birth rate). Both new particles move from the point of splitting independently and according to the same rules as the ones for the initial particle. The life span of each particle up to the reaction (death or splitting) has the law Exp$(\beta+\mu)$.

This model is very close to KPP (Kolmogorov-Petrovskii-Piskunov) model \cite{kpp} with the only difference that the spatial motion (migration of the particles)
in \cite{kpp} is the Brownian motion, and not the random walk with the generator (\ref{HHHa}).

We will focus below on a specific case: $\beta>\mu$ (supercritical process), and the initial population consisting of a single particle located at $x\in R^d$
(it can be a new advanced gene that appeared due to a mutation). Then the first moment (density of population)
\[
m_1(t,x,\Gamma)=E_xn(t,\Gamma)
\]
satisfies the relations
\[
\frac{\partial m_1}{\partial t}=\mathcal Lm_1+(\beta-\mu)m_1,
\]
\[
m_1(0,x,\Gamma)=I_\Gamma(x)=\left\{
                                    \begin{array}{c}
                                      1,~x\in \Gamma ,\\
                                      0,~x\notin\Gamma. \\
                                    \end{array}\right.
\]


Obviously,
\[
m_1(t,x,\Gamma)=\int_\Gamma m_1(t,x,y)dy \quad {\rm and} \quad m_1(t,x,y)=e^{(\beta-\mu)t}p(t,x,y),
\]
where $p(t,x,y)$ is the transition density of the random walk $x(t)$, i.e., $p(t,x,y)dy=P_x\{x(t)\in y+dy\}$, and $p$ is the fundamental solution of the corresponding parabolic problem with the operator $\mathcal L$.

Let us note that functions $m_1$ and $p$ depend on $x-y$, and not on $x,y$ separately, and that they are even with respect to $x-y$ (due to (\ref{aone1sym})). Hence, without  loss of generality, one can assume that $x=0$ (the initial particle is located at the origin). Then we replace variable $y$ by $x$ and introduce the
function $p(t,x)=p(t,0,x)$, the transition density with the initial point at the origin. This function satisfies the relations
\begin{equation}\label{eq1}
\frac{\partial p}{\partial t}=\mathcal L p, \quad p(0,x)=\delta(x).
\end{equation}

The density of the population is large near the origin and decays as $|x|\to \infty$. The later requires some mild assumption on the function $a(z)$. Our results are based on much stronger assumptions on $a(z)$, and the decay of $p(t,x)$ is established below. The front $F_t$ of the propagating population (ecological wave) is defined as the boundary of the set where $m_1\geq 1$:
\begin{equation}\label{front}
F_t=\partial S_t, ~~~ S_t=\{x:~e^{(\beta-\mu) t}p(t,x) \geq 1\}.
\end{equation}

The motion of the front depends on asymptotic behavior of $p(t,x)$ when $t+|x|\to\infty$. This asymptotics ({\it global limit theorem}) is the main goal of the present paper. Asymptotic behavior of $p(t,x)$ depends on the tails of $a(z)$. In the case of very light tails (when $\widehat{a}(k)$ is an entire function of $k$), the picture is similar to the one in the KPP model: the front propagates in time with a constant speed. However, the speed depends on the direction, see \cite{myar}, \cite{kmv}. The proof of the global limit theorem in \cite{myar}, \cite{kmv} is based on a development of the classical Cramer method \cite{IbL}. Random walks on the lattice with heavy tails were studied in \cite{Agmv}, \cite{AsMolV}, 1-D continuous case was studied in \cite{MolPSq}.


This paper concerns the intermediate case of moderate tails. We assume the existence of finitely many moments, including the second one, i.e., $\int_{R^d}|z|^2a(z)dz<\infty$. For the sake of simplicity of formulas, we assume that ${\rm Cov}~ a(\cdot)=I$, i.e.,
\begin{equation}\label{c1}
\widehat{a}(k)=\int_{R^d}e^{-ikz}a(z)dz=1-\frac{|k|^2}{2}+o(|k|^2), \quad |k|\to 0.
\end{equation}
 We will also assume that the intensity of jumps $\chi$ in (\ref{HHHa}) is equal to one, since this can be easily achieved by a time rescaling.
Then the local CLT (Central Limit Theorem) \cite{IbL}, implies that
\begin{equation}\label{clt}
p(t,x)\sim \frac{e^{-\frac{|x|^2}{2t}}}{(2\pi t)^{d/2}}, \quad  |x|<C\sqrt t,~~t \to\infty.
\end{equation}
This statement requires some assumptions on $a(z)$ (additional to the existence of the second moment). They will be automatically valid for our problem. Unfortunately, asymptotics (\ref{clt}) concerns only the ``central zone" $|x|<C\sqrt t$, and we need a global limit theorem that provides asymptotics of $p$ at infinity in the whole $(t,x)$-space. The following example \cite[Ch XIV, s. 6]{IbL} illustrates the desired result in 1-D case.

Let $X_1,...,X_n$ be i.i.d.r.v. (independent and identically distributed random variables) with the  distribution density $a(x)=\frac{2}{\pi(1+x^2)^2},~x\in R.$ Then $EX_i=0$, ${\rm Var} ~X_i$=1. Let $p_n(x)$ be the distribution density of $S_n=X_1+...+X_n $. Then the following relation holds uniformly in $n,|x|$:
\[
p_n(x)=\frac{e^{-\frac{x^2}{2n}}}{(2\pi n)^{1/2}}(1+o(1))+\frac{2nI_{x^2>n}}{\pi(1+x^2)^2}(1+o(1)), \quad n+|x|\to\infty,
\]
where $I_{x^2>n}$ is the indicator of the set $x^2>n$. Under some (very strong) conditions on the tails of $a(\cdot)$, a similar result was proved in \cite{IbL} for more general 1-D walks with moderate tails.

The main requirements in papers \cite{IbL}, \cite{MolPSq} on the 1-D problem, and in papers \cite{Agmv}, \cite{AsMolV} on multidimensional random walks with heavy tails are: high smoothness of $\widehat{a}(k)$ outside of the point $k=0$ and a specific type of singularity at $k=0$. Similar conditions will be imposed in the present paper.

Note that from (\ref{aone1}), (\ref{aone1sym}) it follows that
\begin{equation}\label{c0}
\widehat{a}(0)=1, \quad {\rm Im}~\widehat{a}(k)=0, \quad \widehat{a}(k)\to 0~~{\rm as}~~ |k|\to\infty.
\end{equation}
A natural assumption on the density $a$ with moderate tails would be
\begin{equation}\label{minreg}
a(z)\sim\frac{c_0(\dot z)}{|z|^{d+\alpha}}, \quad |z|\to\infty, \quad \alpha>2,
\end{equation}
where $\dot z=\frac{z}{|z|}$ and $c_0(\dot z)=c_0(-\dot z)$ is a positive continuous function on $S^{d-1}$. Then $a$ has moments of order $[\alpha]$ or $\alpha-1$ if $\alpha$ is non-integer or integer, respectively. In fact, a stronger regularity condition will be
imposed. We assume that $a$ is bounded and has the following asymptotic expansion at infinity:
\begin{equation}\label{ast}
a(z)=\sum_{j=0}^N\frac{c_j(\dot{z})}{|z|^{d+\alpha_j}}+O(\frac{1}{|z|^{2d+\alpha+l}}), \quad z\to\infty, \quad \alpha_0=\alpha<\alpha_1<...<\alpha_N, \quad \alpha>2,
\end{equation}
where $l=1$ if $\alpha>[\alpha], ~l=2$ if $\alpha=[\alpha], ~c_0(\dot z)>0,$ and $c_j(\dot z)$ are sufficiently smooth.

 Expansion (\ref{ast}) immediately implies that  $\widehat{a}(k)$  has the following property:
\begin{equation}\label{c2}
\widehat{a}(k)=\sum_{j=0}^Nb_j(\dot{k}) |k|^{\alpha_j}+a_1(k), ~~a_1(k)\in C^M, \quad M=M(\alpha)=d+[\alpha]+1,
\end{equation}
where $b_j(\dot{k}) |k|^{\alpha_j}$ are the Fourier transforms of $\frac{c_j(\dot{z})}{|z|^{d+\alpha_j}}$.

We also will impose a smoothness condition on $a$. In the lattice version of the problem, it is assumed that
\begin{equation}\label{glad}
b_j(\dot{k})\in C^M, \quad M=d+[\alpha]+1,
\end{equation}
i.e., $b_j$ is $M$ times differentiable on the unit sphere. In the continuous version of the problem (in $R^d$), we require an additional smoothness that guarantees the boundedness and integrability of
 $\widehat{a}(k)$ and its derivatives at infinity:

\begin{equation}\label{c3}
\sup_{|k|\geq 1}|\partial^s_k\widehat{a}(k)|\leq C, \quad \int_{|k|\geq 1}|\partial^s_k \widehat{a}(k)|dk<\infty ,  \quad  |s|\leq M.
\end{equation}
Here $s=(s_1,s_2,...,s_d),~~|s|=\sum s_i,~~\partial^s_k=\frac {\partial^{|s|}}{\partial^{s_1}_{k_1}...\partial^{s_d}_{k_d}}$.

We will use the following property of the characteristic function $\widehat{a}(k)$, which follows immediately from (\ref{aone1}) (since $e^{-ikz}$ is not
identically equal to one when $k\neq 0$):
\begin{equation}\label{cover}
|\widehat{a}(k)|<1 \quad {\rm for} ~0\neq k\in R^d.
\end{equation}
A similar inequality holds also in the lattice case if the random walk is not supported on a sublattice.
\begin{theorem}\label{t}
Let conditions (\ref{aone1sym}), (\ref{ast}), (\ref{glad}), and (\ref{c3}) hold. Then there are constants $A,\varepsilon>0,$ such that the solution of (\ref{eq1}) has the following asymptotic behavior as $ |x|^2\geq At, ~~t\geq 0,~~|x|\to\infty$:
\begin{equation}\label{asu1}
p(t,x)=
     \frac{t}{|x|^{d+\alpha}}[c_0(\dot{x})+O((\frac{1+t}{|x|^2})^\varepsilon)]+E(t,x)(1+O(\frac{t^{1/\alpha}}{|x|})),
    \end{equation}
where $\alpha$ and $c_0(\dot{x})$ are defined in (\ref{ast}), and $E(t,x)=\frac{1}{(2\pi t)^{d/2}}e^{-\frac{|x|^2}{2t}}$.
\end{theorem}
{\bf Remarks.} 1) The first term in (\ref{asu1}) dominates outside of a logarithmic neighborhood of the paraboloid $t=|x|^2$. The second term is larger inside of this neighborhood.

2) We will prove the theorem for the problem in $R^d$, but the statement and the proof remain the same in the lattice case. The following obvious changes are needed to be made for the problem in $Z^d$. Function $\widehat{a}(k)$ will be defined as the Fourier series instead of the Fourier transform. The integration in (\ref{inv1})-(\ref{inv2}), (\ref{derv}), (\ref{ooa}) will be over the torus, not $R^d$, and this makes the assumption (\ref{c3}) on the behavior of $\widehat{a}(k)$ at infinity unnecessary. There is no need for any other changes.

3) Here, and everywhere below, one can replace (\ref{ast}) by a slightly weaker assumption (\ref{c2}).

The following statement is an obvious consequence of Theorem \ref{t}
\begin{theorem}\label{tt}
The front $F_t$ (defined by (\ref{front})) has the form
\[
|x|=(tc_0(\dot{x}))^{\frac {1}{d+\alpha}}e^{\frac {\Delta}{d+\alpha}t}(1+O(t^{-\varepsilon})), \quad t\to\infty.
\]
\end{theorem}

Another application of Theorem \ref{t} concerns the analysis  of the ecological waves and stability of the steady state in the case when $\beta=\mu$ in the presence of a locally supported perturbation: $v(x)=\beta(x)-\mu(x)\in C_0^\infty$, where operator $H$ has the form
\begin{equation}\label{hhhh}
H\psi=\mathcal L\psi+\sigma V(x)\psi, ~~x\in R^d.
\end{equation}
The spectral theory of non-local Schrodinger operators with applications to ecological waves in the case of ultra light tails was developed in   \cite{kmv} (see in particular, Theorem 6.3). If operator $H$ has a single positive eigenvalue $\lambda=\lambda_0(\sigma)$ (with the positive eigenfunction $\psi_0(x)$), then the density of the population in any bounded region grows exponentially with time, but decays at infinity if the population starts with a single particle.  Since the front of the ecological wave is defined by the equation $m_1(t, x)=1$ and
\[
m_1(t, x)\sim e^{\lambda_0t}\psi_0(x), \quad t+|x|\to\infty,
\]
one needs to know the asymptotics of $\psi_0$ at infinity, which coincides with the asymptotics of the Green function
\begin{equation}\label{Gl}
G_\lambda(x)=\int_0^\infty e^{-\lambda t}p(t,x)dt, \quad \lambda>0,
\end{equation}
of the unperturbed operator $\mathcal L$. The latter fact explains the importance of the following statement.
\begin{theorem}\label{ttt}
Let conditions (\ref{aone1sym}), (\ref{ast}), (\ref{glad}), and (\ref{c3}) hold. Then the Green function $G_\lambda(x)$ has the following asymptotic behavior when $0 < \lambda\leq \Lambda_0<\infty,~\lambda |x|^2\to\infty$ (where $\varepsilon$ is defined in (\ref{asu1})):
\begin{equation}\label{Gla}
G_\lambda(x)=
     \frac{1}{\lambda^2|x|^{d+\alpha}}[c_0(\dot{x})+O(\frac{1}{(\lambda|x|^2)^\varepsilon})]+O(\frac{e^
 {-2\sqrt\lambda|x|}}{|x|^{d-2}}).
\end{equation}
In particular, if $\lambda>0$ is fixed, then
\[
G_\lambda(x)=
     \frac{1}{\lambda^2|x|^{d+\alpha}}[c_0(\dot{x})+O(|x|^{-2\varepsilon})], \quad |x|\to\infty.
\]
\end{theorem}
{\bf Remark.} If operator (\ref{hhhh}) has a single positive eigenvalue $\lambda=\lambda_0(\sigma)$, then the front of the ecological wave has the form
\[
|x|\sim(\frac{c_0(\dot{x})}{\lambda_0^2})^{\frac {1}{d+\alpha}}e^{\frac {\lambda_0}{d+\alpha}t}, \quad t\to\infty.
\]
\section{Proofs of the main results}
{\bf Proof of Theorem \ref{t}.}
The solution $p$ of (\ref{eq1}) has the form:
\begin{equation}\label{inv1}
p=\frac{1}{(2\pi)^d}\int_{R^d}e^{ikx+t(\widehat{a}(k)-1)}dk,
\end{equation}
where the integral is understood as the inverse Fourier transform in the sense of distributions, but it also can be reduced to a convergent integral:
\begin{equation}\label{utov}
p=e^{-t}\delta(x)+v(t,x), \quad v=\frac{ e^{-t}}{(2\pi)^d}\int_{R^d}e^{ikx}(e^{t\widehat{a}(k)}-1)dk.
\end{equation}
Indeed, from (\ref{c0}) it follows that $e^{t\widehat{a}(k)}-1\sim ta(k)$ as $k\to\infty$, and the convergence of the integral in (\ref{utov}) follows from (\ref{c3}) with $|s|=0$ and  (\ref{c2}).

First, we are going to justify (\ref{asu1}) when $0\leq t\leq 1$. Let
\begin{equation}\label{hhh}
h=h(t,k)=e^{t(\widehat{a}(k)-1)}-1-t(\widehat{a}(k)-1).
\end{equation}
  Then
\begin{equation}\label{esh}
\int_{R^d}|\partial^j_kh|dk\leq Ct, \quad j=(j_1,...,j_d),~~|j|=M, \quad 0\leq t\leq 1.
\end{equation}
Indeed, each derivative of $h$ contains factor $t$. Estimate (\ref{esh}) for the integral over the region $|k|>1$ follows immediately from (\ref{c3}). The same estimate for the integral over the ball $|k|<1$ is based on (\ref{c2}) and (\ref{glad});  we need only to check that $\partial^j_kh$ has an integrable singularity at $k=0$. From (\ref{c2}) it follows that $\widehat{a}(k)\sim |k|^\alpha, ~k\to 0,$ if one neglects $M$ times differentiable functions. This and  (\ref{c1}) imply that $(\widehat{a}(k)-1)^2\sim |k|^{\alpha+2}, ~k\to 0,$ up to $M$ times differentiable functions. Then the same is true for function $h$. To be more exact, from (\ref{c1}), (\ref{c2}) and (\ref{glad}), it follows that
\[
|\partial^j_kh|\leq t(C_1+C_2|k|^{\alpha+2-|j|}), \quad |k|<1, \quad 0\leq t\leq 1.
\]
This completes the proof of (\ref{esh}) since $|j|=M=d+[\alpha]+1$.

For the inverse Fourier transform $F^{-1}$, we have $|F^{-1}\partial^j_kh|=|x^jF^{-1}h|$. Hence (\ref{esh}) implies that $|F^{-1}h|\leq Ct|x|^{-M}, ~0\leq t\leq 1$. This proves that $p$  is equal to the first term in the right-hand side of (\ref{asu1}) when $0\leq t\leq 1$, since $p=F^{-1}h+(1-t)\delta(x)+ta(x) $ due to (\ref{inv1}) and (\ref{hhh}). The second term in the right-hand side of (\ref{asu1}) is smaller than the remainder in the first term when $0\leq t\leq 1$. Hence, (\ref{asu1})  is proved for $0\leq t\leq 1$.

Our next goal is to prove (\ref{asu1}) for $ t\geq 1$. The reader needs to keep in mind that everywhere below we assume that $t\geq 1$.
The main contribution to the asymptotic behavior of $p$ is given by integration in (\ref{inv1}) over a small neighborhood of the origin. To be more exact, the following statement is valid. Let $\psi=\psi(k)$ be a cut-off function such that $\psi\in C^\infty,~\psi(k)=1$ when $|k|<1/2,~\psi(k)=0$ when $|k|>1$, and let $\psi_1(t,k)=1-\psi(t^{1/\gamma}
k),~2<\gamma<\min(\alpha,3)$. Let
\begin{equation}\label{inv2}
u_1=\frac{1}{(2\pi)^d}\int_{R^d}e^{ikx}e^{t(\widehat{a}(k)-1)}\psi_1(t,k)dk.
\end{equation}
\begin{lemma}\label{l02}
Let conditions (\ref{c1}), (\ref{c2})-(\ref{c3}) hold. Then there exist $\delta>0$ and $C<\infty$ such that
\begin{equation}\label{asu}
|u_1|\leq
    \frac{ C e^{-\delta t^{1-2/\gamma}}}{|x|^{d+[\alpha]+1}},~~ |x|\geq 1,~~ t\geq 1.
    \end{equation}
\end{lemma}
{\bf Proof.} Obviously,
\begin{equation}\label{derv}
(-ix)^ju_1=\frac{1}{(2\pi)^d}\int_{R^d}e^{ikx}\partial^j_k[e^{t(\widehat{a}(k)-1)}\psi_1(t,k)]dk.
 \end{equation}
Thus the lemma will be proved if we show that
\begin{equation}\label{ooa}
\int_{R^d}|\partial^j_k\varphi|dk\leq Ce^{-\delta t^{1-2/\gamma}},~~t\geq 1, \quad  |j|= [\alpha]+d+1, \quad  \varphi:=[e^{t(\widehat{a}(k)-1)}\psi_1(t,k)].
 \end{equation}

Let us estimate the exponential factor in $\varphi$ on the support of $\psi_1$. Relations (\ref{cover}) and (\ref{c0}) imply that the function $\widehat{a}(k)$ is real, vanishes at infinity, and achieves its maximum value $\widehat{a}=1$ at a single point $k=0$. Thus from (\ref{c1}) it follows that there is a constant $\delta>0$ such that $1-\widehat{a}(k)\geq 4\delta|k|^2$ when $|k|\leq 1$. Since $|k|\geq \frac{1}{2}t^{-1/\gamma}$ on the support of $\psi_1$, it follows that $1-\widehat{a}(k)\geq  \delta t^{-2/\gamma}$ on the intersection of the support of $\psi_1$ and the ball $|k|\leq 1$. Hence
\begin{equation}\label{ooo1}
e^{t(\widehat{a}(k)-1)}\leq e^{-\delta t^{1-2/\gamma}}, \quad  k\in \sup \psi_1\bigcap\{k:|k|\leq 1\}, \quad  t\geq 1.
 \end{equation}
Since, $\widehat{a}(k)-1<-\delta_1<0$ for $|k|\geq 1$, we have $e^{t(\widehat{a}(k)-1)}\leq e^{-\delta_1 t}$ when $|k|\geq 1$. From here and (\ref{ooo1}) it follows that
\begin{equation}\label{ooo}
e^{t(\widehat{a}(k)-1)}\leq Ce^{-\delta t^{1-2/\gamma}}, \quad  k\in \sup \psi_1, \quad  t\geq 1.
 \end{equation}

Obviously,
\begin{equation}\label{plpl}
\partial^j_k\varphi=\sum_{|l|\leq|j|}P_l(t,\partial^s_k\widehat{a}(k))e^{t(\widehat{a}(k)-1)}
\partial^l_k\psi_1(t,k),
 \end{equation}
where $P_l$ are polynomials in $t$ and in the derivatives of $\widehat{a}(k)$  with $|s|\leq |j|-|l|$. We also note that
\begin{equation}\label{estp}
|\partial^l_k\psi_1(t,k)|\leq Ct^{|l|/\gamma}.
 \end{equation}

 We split the integral (\ref{ooa}) into two parts: over the ball $|k|<1$ and over the complementary region, and estimate these parts separately. We also take into account that $\varphi=\psi_1=0$ when $|k|t^{1/\gamma}<\frac{1}{2}$. Relations  (\ref{c2}) and (\ref{glad}) imply that
\[
\int_{\frac{1}{2}t^{-1/\gamma}<|k|<1}|P_l(t,\partial^s_k\widehat{a}(k))|dk\leq Ct^m, ~~ t\geq 1,
\]
with some $m<\infty$. From here, (\ref{ooo}), and (\ref{estp}) it follows that
\[
\int_{|k|<1}|\partial^j_k\varphi|dk=\int_{\frac{1}{2}t^{-1/\gamma}<|k|<1}|\partial^j_k\varphi|dk\leq C_1e^{-\delta t^{1-2/\gamma}}, \quad t\geq 1.
\]
A similar estimate for the integral over $|k|>1$ follows from  (\ref{ooo})-(\ref{estp}), (\ref{glad}), (\ref{c3}).  This proves the validity of (\ref{ooa}) and
 completes the proof of the lemma.

\qed

From Lemma \ref{l02} it follows that (\ref{asu1}) needs to be proved only for
\begin{equation}\label{u22}
u_2=\frac{1}{(2\pi)^d}\int_{R^d}e^{ikx}e^{t(\widehat{a}(k)-1)}\psi(t^{1/\gamma}k)dk, \quad t\geq 1.
\end{equation}
Let us note that from (\ref{cover}), (\ref{c0}), (\ref{c1}), and (\ref{c2}) it follows that
\begin{equation}\label{apk0}
\widehat{a}(k)-1=-\frac{|k|^2}{2}+\sum_{j=0}^Nb_j(\dot{k}) |k|^{\alpha_j}+g(k), ~~g(k)\in C^M.
\end{equation}
Here $g(k)=a_1(k)+\frac{|k|^2}{2}$ has zero at $k=0$ of order at least three, i.e.,
\begin{equation}\label{parg}
|\partial^jg(k)|\leq C|k|^{3-|j|}, \quad |k|\leq 1, \quad |j|\leq M.
 \end{equation}
 We add factors $\psi(t^{1/\alpha_j}k/2), ~0\leq j\leq N,$ in the integrand (\ref{u22}). These factors do not change the integrand when $t\geq 1$ since $\gamma<\alpha\leq\alpha_j$,  and therefore $\psi(t^{1/\alpha_j}k/2)=1$ on the support of $\psi(t^{1/\gamma}k)$ when $t\geq 1$. Hence $u_2$ can be rewritten as the following convolution:
\begin{equation}\label{ag}
u_2(t,x)=E*v_0*v_1*...*v_N*v_g,
 \end{equation}
where
\begin{equation}\label{hh}
v_j(t,x)=\frac{1}{(2\pi)^d}\int_{R^d}e^{ikx}e^{tb_j(\dot{k}) |k|^{\alpha_j}}\psi(t^{1/\alpha_j}k/2)dk, \quad v_g(t,x)=\frac{1}{(2\pi)^d}\int_{R^d}e^{ikx}e^{tg(k)}\psi(t^{1/\gamma}k)dk.
 \end{equation}

 Denote by $W_\sigma$ the space of functions $w=w(t,x),~t\geq 1,x\in R^d,$ with the following properties:
\begin{equation}\label{intw}
\int_{R^d}|w|dx<C<\infty, \quad \int_{R^d}w~dx=1,
 \end{equation}
where $C$ does not depend on $t$, and there exist $C_1, ~C_2,~ \varepsilon>0$ such that
 \begin{equation}\label{u222}
w=
    \frac{ t}{|x|^{d+\sigma}} (c(\dot{x})+h (t,x)), \quad  |h(t,x)|<C_1(\frac{t}{|x|^\sigma})^\varepsilon, \quad |x|^\sigma\geq C_2 t, \quad t\geq 1.
    \end{equation}
 \begin{lemma}\label{linc}
 The following inclusions hold: $v_{\alpha_j}\in W_{\alpha_j},~0\leq j\leq N,~~v_g\in W_{\sigma}$ for some $\sigma>\alpha$.
 \end{lemma}
{\bf Proof.} After the substitution $k\to t^{-1/\alpha_j}k $, function $v_j$ takes the form \begin{equation}\label{vjj}
v_j=t^{-d/\alpha_j}q_j(t^{-1/\alpha_j}x),
\end{equation}
where $q_j(x)$ is the inverse Fourier transform of the bounded, compactly supported function $\widehat{q_j}(k)=e^{b_j(\dot{k}) |k|^{\alpha_j}}\psi(k/2)$. Hence, $|q_j(x)|<C.$

We have
\[
\widehat{q_j}(k)=(1+b_j(\dot{k}) |k|^{\alpha_j}+O(b_j^2(\dot{k}) |k|^{2\alpha_j})\psi(k/2).
\]
Since $(1+O(b_j^2(\dot{k}) |k|^{2\alpha_j})\psi(k/2)\in C^M$, the following estimate is valid for the inverse Fourier transform:
\[
|F^{-1}(1+O(b_j^2(\dot{k}) |k|^{2\alpha_j})\psi(k/2)|<C(1+|x|)^{-M}.
\]
Integration $M=d+[\alpha]+1$ times by parts implies that
\[
|F^{-1}(b_j(\dot{k}) |k|^{\alpha_j}(1-\psi(k/2)))|<C|x|^{-M}, \quad |x|\geq1.
\]
Since $F^{-1}(b_j(\dot{k}) |k|^{\alpha_j})=c_j(\dot{x}) |x|^{-d-\alpha_j}$, it follows that
\begin{equation}\label{qjj}
|q_j(x)-c_j(\dot{x}) |x|^{-d-\alpha_j}|<C|x|^{-M}, \quad |x|\geq 1.
\end{equation}
Relations (\ref{qjj}) and (\ref{vjj}) imply (\ref{u222}). The same relations together with the boundedness of $|q_j(x)|$ lead to the first estimate in (\ref {intw}). The second relation in (\ref {intw}) follows from the fact that $v_j$ is the inverse Fourier transform of a function that is equal to one at $k=0$. The statement of the lemma concerning $v_{\alpha_j}$ is proved.

Let us prove that $v_g\in W_{\sigma}.$ The following estimate is valid
\begin{equation}\label{c000}
|\partial^j \widehat{v_g}(t,k)|\leq Ct^{|j|/\gamma}, \quad  |j|\leq M, \quad t\geq 1, \quad \widehat{v_g}(t,k):=e^{tg(k)}\psi(t^{1/\gamma}k).
\end{equation}
In order to prove this estimate, we note that the exponent in the product in (\ref{c000}) is bounded on the support of $\psi$. The latter fact follows from (\ref{parg}) with $|j|=0$ since $\gamma<3$. Hence $| \widehat{v_g}|\leq C$. Next, we note that each differentiation of the product in (\ref{c000}) makes the latter estimate worse at most by factor $Ct^{1/\gamma}$. Indeed, every time when  the derivative is applied to the second factor of the product, the factor $t^{1/\gamma}$ appears. If the derivative is applied to the exponent, then the additional factor $tg'(k)$ appears. Its value on the support of $\psi$ does not exceed $Ct|k|^2\leq Ct^{1-2/\gamma}\leq Ct^{1/\gamma}$ if $t\geq1$. If the derivative is applied to the pre-exponential factor, which was obtained as a result of differentiation of the exponent during the previous steps, then the order with which the pre-exponential factor vanishes at $k=0$ decreases by one. The latter fact also implies the worsening of the estimate for the derivative of $\widehat{v_g}$ by at most $Ct^{1/\gamma}$. Thus (\ref{c000}) is proved.

We estimate $v_g$ in the region $|x|^\gamma<t$ using uniform estimate of the integrand, and we estimate $v_g$ in the region $|x|^\gamma>t$ using integration of the second integral in (\ref{hh}) by parts $M$ times, applying (\ref{c000}), and taking into account the fact that  the  integrand vanishes when $|k|>t^{-1/\gamma}$. This leads to
\[
|v_g|\leq \begin{cases}
     Ct^{-d/\gamma},\quad  |x|^{\gamma}<t,\\
      C\frac{t^{M/\gamma}}{|x|^M}t^{-d/\gamma}, \quad   |x|^\gamma \geq t .
    \end{cases}
\]

This estimate immediately implies the first relation in (\ref{intw}).  The second relation in (\ref{intw}) follows from the fact that $v_g$ is the inverse Fourier transform of a function that is equal to one at $k=0$. In order to obtain (\ref{u222}) for $v_g$, we choose $\sigma$ such that $\alpha<\sigma<[\alpha]+1$ and write it in the form $\sigma=([\alpha]+1)/(1+\varepsilon)$.  Since $\sigma>\alpha>\gamma$, the above estimate on $v_g$ implies that
\begin{equation}\label{egeg}
|v_g|\leq C\frac{t^{([\alpha]+1)/\gamma}}{|x|^{d+[\alpha]+1}}\leq C\frac{t^{([\alpha]+1)/\sigma}}{|x|^{d+[\alpha]+1}}=C\frac{t^{1+\varepsilon}}{|x|^{d+\sigma(1+\varepsilon)}}=
\frac{t}{|x|^{d+\sigma}}O(\frac{t^\varepsilon}{|x|^{\sigma\varepsilon}}), \quad |x|^\gamma\geq t.
\end{equation}
If $|x|^\gamma>t>1$, then $|x|^\sigma>t$, and therefore estimate (\ref{egeg}) justifies the validity of (\ref{u222}) for $v_g$ with $c(\dot{x})=0$, i.e., $v_g\in W_\sigma$.

\qed

\begin{lemma}\label{lcon}
Let $u\in W_\alpha$ with a smooth coefficient $c(\dot{x})$ in (\ref{u222}), and let $v\in W_\beta$ with $\beta>\alpha$. Then $w:=u*v\in W_\alpha$. Moreover, $w$ has the same coefficient $c(\dot{x})$ as function $u$.
 \end{lemma}
{\bf Proof.} For any two functions with well-defined convolution, we have:
\[
\int_{R^d}|u*v|dx\leq \int_{R^d}|u|dx\int_{R^d}|v|dx, \quad  \int_{R^d}u*vdx= \int_{R^d}udx\int_{R^d}vdx.
\]
Hence, the validity of (\ref{intw}) for $u$ and $v$ implies that the same relations hold for $w$. Let us prove the validity of  (\ref{u222}) for $w$. We have
\begin{equation}\label{www}
w=\int_{|y|<|x|/2}u(t,x-y)v(t,y)dy+\int_{|y|>|x|/2}u(t,x-y)v(t,y)dy.
\end{equation}
Let $|x|^\alpha>Ct,~t>1,$ with $C$ being $2^\alpha$ times larger than constants $C_2$ in formulas (\ref{u222}) for $u$ and $v$. Then (\ref{u222}) is valid (even for twice smaller $|x|$) for both $u$ and $v$ since $\beta>\alpha$. We write $u(t,x-y)$ in the first integrand in the right-hand side above  as follows
\begin{equation*}
u(t,x-y)=\frac{ t}{|x-y|^{d+\alpha}} (c(\frac{x-y}{|x-y|})+O(\frac{t}{|x-y|^\alpha})^\varepsilon)=\frac{ t}{|x|^{d+\alpha}}(c(\dot{x})+O(\frac{|y|}{|x|})+O(\frac{t}{|x|^\alpha})^\varepsilon).
\end{equation*}
From (\ref{intw}) and (\ref{u222}) we obtain that
\begin{equation}\label{22}
\int_{|y|<|x|/2}v(t,y)dy=1-\int_{|y|>|x|/2}v(t,y)dy=1-O(\frac{t}{|x|^\beta})=1-O(\frac{t}{|x|^\alpha}),
\end{equation}
and, for large enough $A$,
\[
\int_{|y|<|x|/2}|O(\frac{|y|}{|x|})v(t,y)|dy\leq \frac{C}{|x|} [\int_{|y|<(At)^{1/\beta}}|yv(t,y)|dy+\int_{(At)^{1/\beta}<|y|<|x|/2}|yv(t,y)|dy]
\]
\begin{equation}\label{11}
\leq \frac{C}{|x|}[(At)^{1/\beta}\int_{|y|<(At)^{1/\beta}}|v(t,y)|dy+\int_{|y|>(At)^{1/\beta}}\frac{t}{|y|^{d+\beta-1}}dy]
\leq \frac{C't^{1/\beta}}{|x|}\leq \frac{C't^{1/\alpha}}{|x|}.
\end{equation}
Hence
\[
\int_{|y|<|x|/2}u(t,x-y)v(t,y)dy=\frac{ t}{|x|^{d+\alpha}} (c(\dot{x})+O(\frac{t}{|x|^\alpha})^\varepsilon).
\]
For the second integral in (\ref{www}), we have
\[
|\int_{|y|>|x|/2}u(t,x-y)v(t,y)dy|\leq \sup_{|y|>|x|/2}|v(t,y)|\int_{|y|>|x|/2}|u(t,x-y)|dy
\]
\[
\leq C\sup_{|y|>|x|/2}|v(t,y)|\leq \frac{C't}{|x|^{d+\beta}}.
\]
The last two relations and (\ref{www}) prove (\ref{u222}) for $w$ with $\sigma=\alpha$ and with the same $c(\dot{x})$ as for $u$.

\qed

From Lemmas \ref{linc}, \ref{lcon} it follows that the convolution $w$ of all the factors in (\ref{ag}) except the first one belongs to $W_\alpha$, i.e.,
\begin{equation}\label{lconv}
u_2=E*w, \quad w\in W_\alpha,
\end{equation}
where the coefficient $c(\dot{x})$ in formula (\ref{u222}) for $w$ is equal to $c_0(\dot{x})$ in (\ref{ast}). It remains to show that (\ref{asu1}) holds for the convolution (\ref{lconv}). This can be done using the same type of arguments as in the proof of Lemma \ref{lcon}. We have
\begin{equation}\label{wE}
u_2=\int_{|y|<|x|/2}E(t,x-y)w(t,y)(y)dy+\int_{|y|>|x|/2}E(t,x-y)w(t,y)dy.
\end{equation}
We assume that $|x|^2>At,~t\geq 1,$ with $A$ being $2^\alpha$ times larger than constant $C_2$ in formula (\ref{u222}) for $w$. Then for $E$ in the first integrand above, we have:
\begin{equation*}
E(t,x-y)=\frac{1}{(2\pi t)^{d/2}}e^{-\frac{|x-y|^2}{2t}}=E(t,x)(1+O(\frac{|y|}{|x|})),
\end{equation*}
and the following analogue of (\ref{22}) is valid
\begin{equation*}
\int_{|y|<|x|/2}w(t,y)dy=1-\int_{|y|>|x|/2}w(t,y)dy=1-O(\frac{t}{|x|^\alpha}).
\end{equation*}
For large enough $A$, the following relation,  which is similar to (\ref{11}), is valid:
\[
\int_{|y|<|x|/2}|O(\frac{|y|}{|x|})w(t,y)|dy\leq \frac{C}{|x|} [\int_{|y|<(At)^{1/\alpha}}|yw(t,y)|dy+\int_{(At)^{1/\alpha}<|y|<|x|/2}\frac{t}{|y|^{d+\alpha-1}}dy]
\]
\begin{equation*}
\leq \frac{C}{|x|}[(At)^{1/\alpha}\int_{|y|^\alpha<At}|w(t,y)|dy+
\int_{(At)^{1/\alpha}<|y|<\infty}\frac{t}{|y|^{d+\alpha-1}}dy]
\leq \frac{C't^{1/\alpha}}{|x|}.
\end{equation*}
The last three relations imply that
\begin{equation}\label{wEE}
\int_{|y|<|x|/2}E(t,x-y)w(t,y)dy=E(t,x)(1+O(\frac{t^{1/\alpha}}{|x|})).
\end{equation}

After the substitution $y\to |x|y$, the second integral $I_2$ in (\ref{wE}) takes the form
\[
I_2=\frac{t}{|x|^\alpha(2\pi t)^{d/2}} \int_{|y|>1/2} \frac{ 1}{|y|^{d+\alpha}} (c(\dot{y})+h)e^{-\frac{|x|^2}{2t}|\dot{x}-y|^2}dy, \quad |h|<C(\frac{t}{(|x|y)^\alpha})^\varepsilon.
\]
Obviously, $\nabla_y|\dot{x}-y|^2=0$ only at $y=\dot{x}$. The Laplace method provides the following asymptotic expansion, as $|x|^2/t\to \infty$, of the integral above with $h=0$:
\[
I_2|_{h=0}=\frac{t}{|x|^{d+\alpha}}(c(\dot{x})+O(\frac{\sqrt{t}}{|x|})).
\]
Then we use the Laplace method for $I_2$ when $c(\dot{y})=0$ and $h$ is replaced by its estimate from above. Since $\alpha>2$, this implies:
\[
I_2|_{c(\dot{y})=0}=\frac {Ct}{|x|^{d+\alpha}}(\frac{t}{|x|^2})^\varepsilon.
\]
The last two relations for $I_2$ together with (\ref{wE}) and (\ref{wEE}) justify (\ref{asu1}) for $u_2, ~t\geq 1$.
This completes the proof of the theorem.

 {\bf Proof of Theorem \ref{ttt}.} From (\ref{Gl}), (\ref{asu1}), and (\ref{clt}) it follows that
 \[
 |G_\lambda(x)-\int_0^{|x|^2/A}\frac{t}{|x|^{d+\alpha}}[c_0(\dot{x})+O(\frac{1+t^\varepsilon}{|x|^{2\varepsilon}})] e^{-\lambda t}dt|\leq\int_0^\infty
 \frac{C}{t^{d/2}}e^{-\frac{|x|^2}{2t}}e^{-\lambda t}dt
 \]
\[
\leq \frac{C_1}{|x|^{d-2}}e^
 {-2\sqrt\lambda|x|}, \quad \lambda|x|^2\to\infty.
\]
For the integrals in the left-hand side above, we have, for each $N>0$,
\[
\int_0^{|x|^2/A} te^{-\lambda t}dt=(-\frac{t}{\lambda}-\frac{1}{\lambda^2})e^{-\lambda t}\large|_0^{|x|^2/A}=\frac{1}{\lambda^2}(1+O(\frac{1}{(\lambda|x|^2)^N})), \quad \lambda|x|^2\to\infty,
\]
\[
\int_0^{|x|^2/A} t^{1+\varepsilon}e^{-\lambda t}dt=\frac{1}{\lambda^{2+\varepsilon}}\int_0^{\lambda |x|^2/A} t^{1+\varepsilon}e^{-t}dt
=\frac{1}{\lambda^{2+\varepsilon}}(C+O(\frac{1}{(\lambda|x|^2)^N})), \quad \lambda|x|^2\to\infty,
\]
where $C=\int_0^{\infty} t^{1+\varepsilon}e^{-t}dt$. These three relations imply (\ref{Gla}).

\qed

\end{document}